\documentclass{article}
\usepackage[latin2]{inputenc}
\usepackage{amsmath,amsthm}
\usepackage{amssymb,latexsym}
\usepackage{enumerate}

\usepackage{graphicx}
\begin{document}

\def\dom{{\rm dom}}
\def\proof{\noindent{\sc Proof. }}
\def\qed{\hfill$\Box$}
\def\qed{\hfill\vrule height6pt width6pt depth1pt}
\def\la{\langle}
\def\ra{\rangle}
\def\implies{\longrightarrow}
\def\den{\noindent{Notation.}}
\title{Transfinite inductions producing coanalytic sets}
\author{Zolt\'an Vidny\'anszky\thanks{Partially supported by the Hungarian
Scientific Research Fund grant no. K72655}}

\maketitle
\renewcommand{\thefootnote}{}



\renewcommand{\thefootnote}{\arabic{footnote}}
\setcounter{footnote}{0}

\abstract{A. Miller proved the consistent existence of a coanalytic two-point
set, Hamel basis and MAD family. In these cases the classical 
transfinite induction can be modified to produce a coanalytic set. 
We generalize his result formulating a condition which can be easily applied in such situations. 
We reprove the classical results and as a new application we show that consistently there exists an uncountable coanalytic subset of the plane that intersects every $C^1$ curve in a countable set.}
\newtheorem{theor}{jejj}[section]
\newtheorem{theorem}[theor]{Theorem}
\newtheorem{co}[theor]{Corollary}
\newtheorem{all}[theor]{Proposition}
\newtheorem{lemma}[theor]{Lemma}
\newtheorem{krd}[theor]{Question}
\newtheorem{example}[theor]{Example}
\newtheorem{defin}[theor]{Definition}
\newtheorem{pr}[theor]{Problem}
\newtheorem{claim}{Claim}
\newtheorem{remark}[theor]{Remark}
\newcommand{\thm}[2]{\begin{theorem}\label{#1}#2\end{theorem}}
\newcommand{\cor}[2]{\begin{co}\label{#1}#2\end{co}}
\newcommand{\prop}[2]{\begin{all}\label{#1}#2\end{all}}
\newcommand{\lem}[2]{\begin{lemma}\label{#1}#2\end{lemma}}
\newcommand{\prob}[2]{\begin{pr}\label{#1}#2\end{pr}}
\newcommand{\ex}[2]{\begin{example}\label{#1}#2\end{example}}
\newcommand{\defi}[2]{\begin{defin}\label{#1}#2\end{defin}}
\newcommand{\rem}[2]{\begin{remark}\label{#1}#2\end{remark}}
\newcommand{\que}[2]{\begin{krd}\label{#1}#2\end{krd}}
\newcommand{\prb}[2]{\begin{pr}\label{#1}#2\end{pr}}
\newcommand{\clm}[2]{\begin{claim}\label{#1}#2\end{claim}}
\newcommand{\R}{\mathbb{R}}
\newcommand{\md}{\models}
\newcommand{\om}{\omega^\omega}
\section{Introduction}
A two-point set is a subset of the plane that intersects every line in exactly two points. Mazurkiewicz showed the existence of a two-point set using transfinite induction.
Erd\H os asked whether a two-point set can be a Borel set. This question is still open.

 A. Miller proved in \cite{miller1} that under certain set theoretic assumptions
(namely $V=L$, where $L$ denotes G\"odel's constructible universe) one can
construct a coanalytic two-point set. Miller also proved the
consistent existence of a coanalytic MAD family and Hamel basis. The author
proves the statement solely for two-point sets and the proof uses deep set
theoretical tools. References to Miller's method appear in several papers
(\cite{vera}, \cite{gao}, \cite{bart} etc.), sometimes omitting the proof. However, the first version of the method was published by Erd\H os, Kunen and Mauldin 
(\cite{erdos}).

Our
aim here is to make precise and prove a "black box" condition which could
easily be applied without the set theoretical machinery.

Let us remark here that in all of the above mentioned cases, except of course
the two-point set, the class of coanalytic sets is best possible, since it
is known that there is no analytic

\begin{enumerate}
 \item MAD family,
  \item Hamel basis,
 \item $C^1$-small set (that is, an uncountable subset of the plain that
intersects every $C^1$ curve in countably many points).
\end{enumerate}

1. is a classical result of Mathias (\cite{mat}) and for the proof of 3.
see \cite{kun}. 2. can be shown with an easy computation. Moreover, assuming
projective determinacy one can show that there is no projective Hamel basis or
$C^1$-small set. It is also an interesting fact that an analytic two-point set
is
automatically Borel.  

Now to formulate our results first we define Turing reducibility. Throughout
the paper $M$ will stand for $\R^n$, $2^{\omega}$,
$\mathcal{P}(\omega)$ or $\omega^{\omega}$.
\defi{}{Suppose that $x,y \in M$. We say that $x$ is Turing reducible to $y$ if there exists a Turing machine that computes $x$ with the oracle $y$. This relation is denoted by $x \leq_T y$. Let us say that $A \subset M$ is cofinal in the Turing degrees, if for every $x \in M$ there exists a $y \in A$ such that $x \leq_T y$.}

Roughly speaking, the theorem will state that if given a transfinite induction
that picks a real $x_{\alpha}$ at each step $\alpha$, the set of possible
choices (described by the set $F$ below) is nice enough and cofinal in the
Turing degrees then the induction can be realized so that it produces a
coanalytic set. In most cases there will be an extra requirement that
$x_{\alpha}$ has to be picked from a given set $H_{\alpha}$. For example, in the
construction of the two-point set $H_{\alpha}$ is the $\alpha^{th}$ line.
Instead of the sets $H_{\alpha}$ we will use a parametrization where
$H_{\alpha}$ will be coded by $p_{\alpha}$ and typically the codes will range
over $\R$. The set of the codes will be denoted by $B$.

\textit{Notation.} If $S \subset X \times Y$ and $x \in X$ we denote the $x$-section of $S$ (i. e. $\{y \in Y: (x,y) \in S\}$) with $S_x$. Let $\omega$ denote the first infinite ordinal, $
\omega_1$ is the first uncountable ordinal. For a set $H$ the set of countable
sequences of elements of $H$ is denoted by $H^{\leq \omega}$. Note that if $M$
is a Polish space then there is a natural Polish structure on $M^{\leq \omega}$.

\defi{}{Let $F \subset M^{\leq \omega} \times B \times{M}$, and $X \subset M$.
We say that $X$ is compatible with $F$ if there exist enumerations
$B=\{p_{\alpha}:\alpha<\omega_1\}$, $X=\{x_{\alpha}:\alpha<\omega_1\}$ and for
every $\alpha<\omega_1$ a sequence $A_{\alpha} \in M^{\leq \omega}$ that is an
enumeration of $\{x_{\beta}:\beta<\alpha\}$ in type $\leq \omega$ such that
$(\forall \alpha<\omega_1)(x_{\alpha} \in F_{(A_{\alpha},p_{\alpha})})$ holds. }

This definition is basically describing that in each step of the transfinite induction we pick an element from a set $F_{(A_{\alpha},p_{\alpha})}$ which depends on the set of the previous choices $A_{\alpha}$ and the $\alpha^{th}$ parameter $p_{\alpha}$. 

\thm{alap}{($V=L$) Let $B$ be an uncountable Borel subset of an arbitrary
Polish space. Suppose  that $F \subset M^{\leq \omega} \times B
\times M$ is a coanalytic set and for all $p \in B$, $A \in M^{\leq \omega}$
the section $F_{(A,p)}$ is cofinal in the Turing degrees.
Then there exists a coanalytic set $X$ that is compatible with $F$.}

In fact we will prove a much stronger theorem (Theorem \ref{main}), which we
call the Main Theorem. However, all the classical applications are using Theorem
\ref{alap} and it will be an easy consequence of the Main Theorem (see Section
4). We would like to emphasize one of our further results from Section 4.

\thm{UDC0}{($V=L$) Suppose that $G \subset \R \times \R^n$ is a Borel set and
for every countable $A \subset \R$ the complement of the set $\cup_{p \in A}
G_p$ is cofinal in the Turing degrees. Then there exists an uncountable
coanalytic set $X \subset \R^n$ that intersects for every $p \in \R$ the section
$G_p$ in a countable set.}                                            

Our paper is organized as follows: in Section 2 we summarize the most important
facts used for the proof and Section 3 contains the proof of the Main Theorem.
In Section 4 we prove several generalizations, a partial converse and we
obtain the existence of a coanalytic Hamel basis (which slightly differs from
the other applications). Finally in Section 5 we present the applications of our
theorem and mention some open problems. The reader only interested in how to
apply the method developed in this paper
may now proceed to Section 5 which is not building on Sections 2, 3 and 4.

\section{Preliminaries}
We will use standard notation as in ~\cite{mv}. If $A$ is a set,
$\mathcal{P}(A)$ denotes the power set of $A$. We identify $\omega^{\omega}$,
$(\omega^{\omega})^{\leq  \omega}$, $2^\omega$,$\omega^\omega$, $\R^{\leq
\omega}$, $\mathcal{P}(\omega)$ and their finite products, since there are
recursive Borel-isomorphisms between them (\cite[3I.4.Theorem]{mv}). A ``real``
is an element of one of
these spaces. For convenience we will use $\omega^\omega$ in most cases. If $A
\in (\om)^{\leq \omega}$ and $n \in \omega$, let us denote the $n^{th}$ element
of $A$ (as a sequence) with $A(n)$. 

As usual, the continuous images of Borel sets are called analytic sets and their
complements are called coanalytic sets. If $t$ is a real, let us denote the
classes of the arithmetic and projective hierarchy recursive in $t$ with
$\Sigma^i_j(t)$, $\Delta^i_j(t)$ and $\Pi^i_j(t)$ ($i=0,1$, $j \in \omega$).
Thus for example the set of coanalytic subsets of $\om$ equals to $\bigcup_{t
\in \om} \Pi^1_1(t)$. For $t=\emptyset$ we will write $\Sigma^i_j$ instead of
$\Sigma^i_j(t)$ etc.

The theorems we will use can be found in \cite{sacks} and \cite{devlin}, but we
recall the most important facts. Let us denote the set of self-constructible
reals, i.e. $\{x \in \omega^\omega: x \in L_{\omega^x_1}\}$ with $\mathcal{S}$,
where $\omega^x_1$ is the first ordinal not recursive in $x$ and $L_\alpha$ is
the $\alpha^{th}$ level of G\"odel's constructible universe, $L$. Let $<_L$ be
the standard well ordering of $L$.

\thm{a}{(\cite[Theorem (2A-1)]{kechris}) $\mathcal{S}$ is a $\Pi^1_1$ set.}
For reals $x,y$ let us denote by $x \leq_h y$ that $x$ is hyperarithmetic in
$y$ or equivalently $x \in \Delta^1_1(y)$ (see \cite{sacks} or \cite[Corollary 27.4]{millerk}
). If $A$ is a set, $L_{\alpha}[A]$ denotes the $\alpha^{th}$
level of the universe constructed from $A$, that is, in the initial step we
start from $\emptyset$ and $A$. 
\thm{b}{(\cite[Part A, Chapter II, 7]{sacks}) $x \leq_h y$ is a $\Pi^1_1$
relation and for arbitrary reals it is equivalent to $x \in L_{\omega^y_1}[y]$.
Moreover, $x \leq_h y$ implies $\omega^x_1 \leq \omega^y_1$.}
We will use the following form of Spector-Gandy-theorem:
\thm{sck}{(\cite[Corollary 29.3]{millerk}) Let $A \subset (\omega^\omega)^2$ be a $\Pi^1_1(t)$ subset of $(\omega^\omega)^2$. Then the set \[(\exists y \leq_h x)((x,y) \in A)\]
is also $\Pi^1_1(t)$.}
In \cite{kin} the authors work with a very useful alternative form. We call a formula in the language of set theory $\Sigma_1$ if it has just one unbounded quantifier and that is existential. In case all the quantifiers are bounded, we call it $\Delta_0$.  
\thm{spg}{A set $A$ is $\Pi^1_1(t)$ if and only if there exists a $\Sigma_1$ formula $\theta$ such that
\[ x \in A \iff L_{\omega^{(x,t)}_1}[x,t]\models \theta(x,t).\]}

\defi{cofprogd}{We call a set $X \subset \omega^\omega$ cofinal in the hyperdegrees if for every $y \in \omega^\omega$ there exists an $x \in X$ such that
$y \leq_h x$.}
Furthermore, in \cite{kin} one can find the following lemma.
\lem{cofprog}{($V=L$) Let $t \in \omega^\omega$ be arbitrary. A $\Pi^1_1(t)$
set $X$ is cofinal in the hyperdegrees if and only if $X \cap \mathcal{S}$ is
cofinal in the hyperdegrees.}

\section{The main theorem}

First we will prove a rather technical lemma.
\lem{TEK}{Suppose that $\theta(s,p,q)$ is a $\Sigma_1$ formula of set theory. Then there exists a $\Sigma_1$ formula $\theta'(s,p)$ such that for every limit ordinal $\alpha>\omega$ \[L_{\alpha} \models ((\forall q <_L p)(\theta(s,p,q)) \iff \theta'(s,p)).\]}

\proof By \cite[3.5 Lemma, p. 75]{devlin}  there exists a $\Sigma_1$ formula
$\zeta(x,y)$ such that for arbitrary limit ordinal $\alpha>\omega$ and $x,y
\in L_{\alpha}$
\[L_{\alpha} \models (\zeta(x,y) \iff y=\{t:t <_L x\}).\]
Notice that if $\alpha>\omega$ is a limit ordinal and $x \in L_{\alpha}$ then $\{t:t <_L x\} \in L_{\alpha}$.
Let \[ \theta''(s,p)=(\exists y)(\zeta(p,y) \land (\forall q \in y)(\theta(s,p,q))).\] Now, since $\theta''$ contains solely existential and bounded quantifiers, using the well-known trick there exists a $\Sigma_1$ formula $\theta'(s,p)$ such that for every limit ordinal $\alpha>\omega$ \[L_{\alpha} \models (\theta''(s,p) \iff \theta'(s,p)).\]
\qed

In the following lemma we will select a single well-ordering of $\omega$ of type
$\alpha$ for every countable ordinal $\alpha$ in a "nice" way. The selection
will be done by a formula $\phi(z,x)$ that intuitively means that $x$ "knows"
that $z$ is a canonical well-ordering. 
Let $z \subset \omega^2$ and define $<_z$ as the relation $m <_z n \iff (m,n) \in z$.  Let us use the notation $dom(<_z)$ for the set $\{n \in \omega: (\exists m \in \omega)((m,n) \in z)\}$. For $z,z' \in \mathcal{P}(\omega^2)$ we say that $<_z \cong <_{z'}$ if there exists a bijection $f:dom(<_z) \to dom(<_{z'})$ such that \[(\forall m, n \in dom(<_z))(m <_z n \iff f(m)<_{z'} f(n)).\]
Now if $<_z$ is an ordering and $n \in \omega$ let us denote by $<_z|_{<_z n}$ the ordering obtained by restricting $<_z$ to the set $\{m \in \omega:m <_z n\}$.
\lem{fi}{($V=L$) There exists a formula $\phi(z,x)$ defining a $\Pi^1_1$ subset
of $\mathcal{P}(\omega^2) \times \om$ with the following properties
\begin{enumerate}
 \item if $s \subset \omega^2$ and $<_{s}$ is a well-ordering then there exists a unique $z$ such that $<_z \cong <_{s}$, $(\exists x \in \om) \phi(z,x)$ and $dom(<_z)$ is a natural number or $\omega$ 
\item if $y \in \mathcal{S}$, $x \leq_h y$ and $\phi(z,x)$ then $\phi(z,y)$
\item if $\phi(z,x)$ then $z \leq_h x$ and $x \in \mathcal{S}$
\item if $\phi(z,x)$ and $n \in \omega$ is arbitrary then there exists a unique pair $g_n,y_n \in L_{\omega^x_1}$ such that $\phi(y_n,x)$ and $g_n \subset \omega^2$ is an isomorphism between $<_z|_{<_z n}$ and $<_{y_n}$. 

\end{enumerate}}
\proof
First let us denote by $\psi(z,h,\alpha)$ the conjunction of the following three
formulas:
\begin{itemize}
\item $h$ is a function, $dom(h)=\alpha$ is an ordinal, $ran(h)=dom(<_z)$
\item $(\forall \beta,\beta' \in \alpha) (\beta \in \beta' \iff h(\beta)<_z h(\beta'))$
\item $dom(<_z)$ is a natural number or $\omega$.
\end{itemize}
 So $\psi(z,h,\alpha)$ says that $h$ is an isomorphism between $\alpha$ and
$<_z$. Notice that $\psi$ is a $\Delta_0$ formula (see \cite{devlin}, Section
I). Hence for limit ordinals $\beta>\omega$ if $z,h,\alpha \in
L_{\beta}$ then
$L \models \psi(z,h,\alpha) \iff L_{\beta} \models \psi(z,h,\alpha)$.

Let us define $\phi(z,x)$ as follows: \[\phi(z,x)  \iff x \in \mathcal{S} \land
z \leq_h x \land \] \[L_{\omega^x_1} \models (\exists h \exists
\alpha)\big((\psi(z,h,\alpha)  \land (\forall (z',h') <_L (z,h))(\lnot
\psi(z',h',\alpha))\big).\]

First, we will prove that $\phi(z,x)$ defines a $\Pi^1_1$ set. 
The formula \[(\exists h \exists \alpha)\big((\psi(z,h,\alpha)  \land (\forall (z',h') <_L (z,h))(\lnot \psi(z',h',\alpha))\big)\]
by Lemma \ref{TEK} is equivalent to a $\Sigma_1$ formula, say $\zeta(z)$, in $L_{\beta}$ if $\beta$ is a limit ordinal and $\beta>\omega$. 
Notice that $z \leq_h x$ implies $(x,z) \leq_h x$ so $\omega^{(x,z)}_1 \leq
\omega^x_1$ by Theorem \ref{b}. Moreover, from $(x,z) \leq_h x$ and by Theorem
\ref{b} we have that $(x,z) \in L_{\omega^x_1}[x]$. Additionally, $x \in
\mathcal{S}$ so $L_{\omega^x_1}=L_{\omega^x_1}[x]$. Thus $(x,z) \in
L_{\omega^x_1}$ and the equality $L_{\omega^{(x,z)}_1}[x,z]=L_{\omega^x_1}$
holds. Therefore 
\[L_{\omega^x_1} \models (\exists h \exists \alpha)\big((\psi(z,h,\alpha)  \land (\forall (z',h') <_L (z,h))(\lnot \psi(z',h',\alpha))\big)\]
\[\iff\]
\[L_{\omega^x_1} \models \zeta(z)\]
\[\iff\]
\[L_{\omega^{(x,z)}_1}[x,z] \models \zeta(z).\]

By Theorems \ref{a} and \ref{b} it is clear that $(x \in \mathcal{S}) \land (z
\leq_h x)$ defines a $\Pi^1_1$ set. Now we can prove that the set
$\{(x,z):L_{\omega^{(x,z)}_1}[x,z] \models \zeta(z)\}$ is also $\Pi^1_1$ using
Theorem \ref{spg} with $t=0$ and replacing $x$ by $(x,z)$. Thus $\phi$ defines a
$\Pi^1_1$ set. 

Now we will prove that $\phi(z,x)$ has the required properties. 
\begin{enumerate}
\item Let $s \subset \omega^2$ be an arbitrary well-ordering.
Then $<_{s}$ is isomorphic to some ordinal $\alpha$. There exists a $<_L$ minimal pair $(z,h)$ such that $h$ is an isomorphism between $<_z$ and $\alpha$ and $\dom(<_z)$ is a natural number or $\omega$. Therefore 
\[L \models (\exists h \exists \alpha)\big((\psi(z,h,\alpha)  \land (\forall (z',h') <_L (z,h))(\lnot \psi(z',h',\alpha))\big).\]
Notice that if $\xi(s)$ is a $\Delta_0$ formula, $\beta$ is a limit ordinal such
that $s \in L_{\beta}$ and $L \models \xi(s)$ then $L_{\beta} 
\models \xi(s)$. Therefore automatically $L_{\beta} \models (\exists r)(\xi(r))$. Considering this one can conclude that \[L_{\omega^x_1} \models (\exists h \exists \alpha)\big((\psi(z,h,\alpha) \land (\forall (z',h')<_L (z,h))(\lnot \psi(z',h',\alpha))\big)\] 
holds if $(z,h) \in L_{\omega^x_1}$. $\mathcal{S}$ is cofinal in the hyperdegrees (Lemma \ref{cofprog}) hence 
there 
exists an $x \in \mathcal{S}$ such that $(z,h) \in L_{\omega^x_1}$. So for such an $x$ we have $\phi(z,x)$.

\item To prove the second claim just observe that 
$\Sigma_1$ formulas are upward absolute for transitive sets and notice that $x \leq_h y$ implies that $L_{\omega^x_1} \subset L_ {\omega^y_1}$.
\item Obvious from the definition of $\phi$.

\item Let $x \in \om$, $z \subset \omega^2$, $n \in \omega$ and assume that $\phi(z,x)$ holds. Clearly there exists a unique ordinal $\beta<\alpha$ such that $\beta \cong <_z|_{<_z n}$. 

First we will prove that there exists a pair $(y'_n,h'_n) \in \mathcal{P}(\omega^2) \times \beta^{dom(<_{y'_n})}$ so that $L_{\omega^x_1} \models \psi(y'_n,h'_n,\beta)$. We know that $L_{\omega^x_1} \models 
\psi(z,h,\alpha)$ for some $h, \alpha \in L_{\omega^x_1}$ so the same holds in $L$. The fact that $\psi(z,h,\alpha)$ holds implies that $h$ is an isomorphism between $<_z$ and $\alpha$, so $h'=h|_{\beta}$ is an isomorphism between $\beta$ and $<_z|_{<_z n}$. Obviously, $h' \in L_{\omega^x_1}$, so there exists an ordinal $\gamma<\omega^x_1$ such that $h' \in L_{\gamma}$. 

Let $e: \omega \to ran(h')$ be defined as follows: 
\[\langle m,k  \rangle \in e \iff (k \in ran(h') \land \exists e'(e':m \leftrightarrow ran(h') \cap k+1)),\]  
in other words, there exists a bijection between $m$ and the initial segment of $ran(h')$, or equivalently, $|\{l \in ran(h'):l \leq k\}|=m$.
Since the bijections between the finite subsets of $\omega$ are already in $L_{\omega}$, we have that $e \in L_{\gamma+2} \subset L_{\omega^x_1}$. $e$ is clearly a one-to-one function from a finite number or $\omega$ onto $ran(h')$. 

Now take $\langle k,l \rangle \in y'_n \iff \langle e(k),e(l) \rangle \in z$ and $h'_n=e^{-1} \circ h'$. Then $L \models \psi(y'_n,h'_n,\beta)$ and of course $y'_n, h'_n, \beta \in L_{\omega^x_1}$ hence $L_{\omega^x_1} \models \psi(y'_n,h'_n,\beta)$. 

Thus there exists a $<_L$ minimal pair $(y_n,h_n) \in L_{\omega^x_1}$ such that $L_{\omega^x_1} \models \psi(y_n,h_n,\beta)$. Note that the $<_L$ ordering is absolute for $L_{\alpha}$ and $L$ if $\alpha>\omega$ is a limit ordinal, so $L_{\omega^x_1} \models "(y_n,h_n)$ is the $<_L$ minimal pair such that $\psi(y_n,h_n,\beta)"$. By Theorem \ref{b}, if $y_n \in 
L_{\omega^x_1}$ then $y_n \leq_h x$. Thus $\phi(y_n,x)$ holds.  

Finally recall that $h_n:\beta \to dom(<_{y_n})$ and $h':\beta \to dom(<_z|_{<_z n})$ are isomorphisms in $L_{\omega^x_1}$. So the function $g_n=h_n \circ (h')^{-1}$ is in $L_{\omega^x_1}$. This is an isomorphism between two well-orderings so this is unique. 
\end{enumerate}
\qed
\\
\\
Let us recall the definition of compatibility.

\defi{}{Let $F \subset M^{\leq \omega} \times B \times{M}$, $X \subset M$. We
say that $X$ is compatible with $F$ if there exist enumerations
$B=\{p_{\alpha}:\alpha<\omega_1\}$, $X=\{x_{\alpha}:\alpha<\omega_1\}$ and for
every $\alpha<\omega_1$ a sequence $A_{\alpha} \in M^{\leq \omega}$ that is an
enumeration of $\{x_{\beta}:\beta<\alpha\}$ in type $\leq \omega$ such that
$(\forall \alpha<\omega_1)(x_{\alpha} \in F_{(A_{\alpha},p_{\alpha})})$ holds. }

\thm{main}{(Main Theorem) ($V=L$) Let $t \in \om$. Suppose that $F \subset
(\omega^\omega)^{\leq \omega} \times \omega^\omega \times \omega^\omega$ is a
$\Pi^1_1(t)$ set and for all $p \in \omega^\omega$, $A \in (\omega^\omega)^{\leq
\omega}$ the section $F_{(A,p)}$ is cofinal in the hyperdegrees.
Then there exists a $\Pi^1_1(t)$ set $X \subset \omega^{\omega}$ that is compatible with $F$.}
\noindent{\sc Proof of the Main Theorem.} 

In the first step we will modify the set $F$.
Let us define \[F' \subset \mathcal{P}(\omega^2) \times (\om)^{\leq \omega} \times (\om)^{\leq \omega} \times \om \times \om,
(z,A,P,p,x) \in F' \iff \]
\begin{enumerate}
\item $\phi(z,x)$ (in particular $x \in \mathcal{S}$)
\item $A,P,p,t \leq_h x$, $(A,p,x) \in F$
\item $L_{\omega^x_1} \models \exists g$
\begin{enumerate}
\item $g$ is a function, $dom(g) \in \omega \cup \{\omega\}$, $ran(g)=P$
\item $(\forall n, m \in dom(g))$ $(n <_z m  \iff g(n) <_L g(m) )$
\item $(\forall p' <_L p)(p' \in \omega^\omega \Rightarrow (\exists n \in \omega)(g(n)=p'))$
\end{enumerate}
\end{enumerate}

The role of $z$ is that it will encode the history of the previous choices. $1-2$ basically ensures that $x$ is complicated enough. The clauses $(a)$ and $(b)$ describe that $P$ is an enumeration in type $\leq \omega$ of the first $\alpha$ reals with respect to $<_L$ where $\alpha = tp(<_z)$. $(c)$ is the formalization of $L_{\omega^x_1}\models$ "$p$ is the $\alpha^{th}$ real with respect to $<_L$".

Lemma \ref{fi}, Theorems \ref{a} and \ref{b} guarantee that the $1$ and $2$ are
defining a $\Pi^1_1(t)$ set. 

We can prove that $3$ defines a $\Pi^1_1$ set similarly as we did in Lemma \ref{fi}: $(a)$ and $(b)$ are $\Delta_0$ formulas, $(c)$ is $\Sigma_1$ by Lemma \ref{TEK}. So by the well-known technical trick the conjunction is equivalent to a $\Sigma_1$ formula. Moreover we know that for arbitrary reals $a \leq_h b \iff a \in L_{\omega^b_1}[b]$ and $a \leq_h b $ implies $\omega^a_1\leq \omega^b_1$. Therefore by $1$ and $2$ \[L_{\omega^{(z,A,P,p,t,x)}_1}[z,A,P,p,t,x]=L_{\omega^x_1}\] and using the Spector-Gandy Theorem (Theorem \ref{spg}) we can conclude that $F'$ is a $\Pi^1_1(t)$ set.

\rem{meg}{By absoluteness, if $(z,A,P,p,x) \in F'$  then $P$
must be
the enumeration of the first $\alpha$ reals given by $<_z$ in $L$ as well.
Similarly $p$
must be the $\alpha^{th}$ real with respect to $<_L$ (where $\alpha=tp(<_z)$).
}

\lem{l1}{Suppose that $x \in F'_{(z,A,P,p)}$, $x \leq_h y$ and $y \in \mathcal{S} \cap F_{(A,p)}$. Then $y \in F'_{(z,A,P,p)}$.}
\proof Let $x, y$ be reals satisfying the conditions above. Now considering the
definition of $F'$, the formula $\phi(z,y)$ holds by the second claim of Lemma
\ref{fi}. Of course, $A,P,p,t \leq_h x$ implies $A,P,p,t \leq_h y$. Finally,
$L_{\omega^x_1} \subset L_{\omega^y_1}$, by Theorem \ref{b}, and the formula
in 3. that
must hold in $L_{\omega^y_1}$ does not depend on $x$, hence it is also true in
$L_{\omega^y_1}$. \qed

\lem{l2}{If the section $F'_{(z,A,P,p)}$ is non-empty then it is cofinal in the hyperdegrees.}
\proof Fix an arbitrary $s \in \om$ and let $x \in F'_{(z,A,P,p)}$. By the assumptions of the Main Theorem each section $F_{(A,p)}$ cofinal in the hyperdegrees.  Using Lemma \ref{cofprog} we have that there exists a $y \in F_{(A,p)} \cap \mathcal{S}$ such that $s,x \leq_h y$. Thus by the previous lemma $y \in F'_{(z,A,P,p)}$ and this proves the statement. \qed

Now we select a real from each nonempty section of $F'$. Let $F'' \subset F'$ be
a $\Pi^1_1(t)$ uniformization of $F'$, that is, for all  $(z,A,P,p) \in proj(F')$ we have
$|F''_{(z,A,P,p)}|=1$ (see \cite{millerk} or \cite{sacks} for the relative
version of the uniformization theorem). 

There may be elements $(z,A,P,p,x) \in F''$ with "wrong" history, namely $A(n)$ may not be a selected real for some $n \in \omega$. So we have to sort out the appropriate ones. 

Let $F''' \subset F''$ be defined as follows:

$(z,A,P,p,x) \in F''' \iff $
\begin{enumerate}
\item $(z,A,P,p,x) \in F''$
\item $(\forall n \in \omega)(\exists g_n,y_n \leq_h x)$
\begin{enumerate}
\item $\phi(y_n,x)$
\item $g_n$ is an isomorphism between $<_z|_{<_z n}$ and $y_n$
\item if $A_n,P_n \in (\om)^{\leq \omega}$ is defined by $A_n(i)=A(g_n(i))$ and similarly $P_n(i)=P(g_n(i))$ then $(y_n,A_n,P_n,P(n),A(n)) \in F''$
\end{enumerate}
\end{enumerate}

By properties of $\phi$, for every countable ordinal $\alpha$ we have a canonical enumeration of $\alpha$. In the definition above (c) ensures that for every $(z,x,A,P,p) \in F'''$ the set $A$ is the canonical enumeration of the previous choices given by the uniformization of $F'$.

The clauses $(a),(b)$ are defining a $\Pi^1_1(t)$ set. Now take the map $\Psi: (A,P,y_n,g_n,n) \mapsto (y_n,A \circ g_n,P \circ g_n,P(n),A(n))$. Observe that $\langle (A,P,y_n,g_n,n), (w_1,w_2,w_3,w_4,w_5) \rangle  \in \Psi \iff $
$y_n=w_1$, $w_4=P(n)$, $w_5=A(n)$ and $(\forall m \in \omega)(w_2(m)=A(g_n(m))
\land w_3(m)=P(g_n(m))$. So $\Psi$ is a $\Delta^1_1$ map and condition $(c)$
describes that $(A,P,y_n,g_n,n) \in \Psi^{-1}(F'')$ thus defines a $\Pi^1_1(t)$
set. Therefore, using Theorem \ref{sck} we can conclude that $F'''$ is also a
$\Pi^1_1(t)$ set. 

Now we will prove that $F'''$ contains a "good selection" and then $X$ will be the projection of $F'''$ on the last coordinate.

More precisely, let: \[x \in X \iff (\exists (z,A,P,p) \leq_h x) ((z,A,P,p,x) \in F''').\]

Notice that $X$ is indeed the projection of $F'''$ on the last coordinate: if $(z,A,P,p,x) \in F''' \subset F'$ then $(A,P,p) \leq_h x$ by the definition of $F'$ and from the $3^{rd}$ point of Lemma \ref{fi} we obtain that $z \leq_h x$, so obviously $ (z,A,P,p) \leq_h x$ holds. 

Observe that by Theorem \ref{sck} the set $X$ is also $\Pi^1_1(t)$. 

\prop{aaa}{For every $\alpha \in \omega_1$ there exists a unique
$(z_{\alpha},A_{\alpha},P_{\alpha},p_{\alpha},x_{\alpha}) \in F'''$ such that
$<_{z_{\alpha}} \cong \alpha$. Moreover, $\{A_{\alpha}(n):n
\in
\omega\}=\{x_{\beta}:\beta<\alpha\}$ holds for every $\alpha<\omega_1$.}

\noindent{\sc Uniqueness.} Let $(z,A,P,p,x),(z',A',P',p',x')\in F'''$ be such that $<_z \cong <_{z'} \cong \alpha$. 

$z=z'$: follows form the $1^{st}$ point of Lemma \ref{fi} since both of $\phi(z,x)$ and $\phi(z',x')$ must hold. 

$p=p'$: clear by Remark \ref{meg}.

$P=P'$: also from Remark \ref{meg} we have that $P$ and $P'$ are
enumerations of the first $\alpha$ reals given by $<_z=<_{z'}$.

$A=A'$: suppose
not. Then take the $<_z$ minimal $n \in \omega$ such that $A(n) \not =A'(n)$. By
the definition of $F'''$ there exist $y_n, g_n$ and $y'_n,g'_n$ 
such that  $(y_n,A_n,P_n,P(n),A(n)) \in F''$ and $(y'_n,A'_n,P'_n,P'(n),A'(n))
\in F''$, $g_n$ and $g'_n$ are isomorphism between $<_z|_{<_z n}$ and $y_n,
y'_n$ 
and $\phi(y_n,x)$ and $\phi(y'_n,x)$ hold. Then again by Lemma \ref{fi} $y_n=y'_n$, $g_n$ is unique so it must be equal to $g'_n$. We obtain that 
$(y_n,A_n,P_n,P(n))=(y'_n,A'_n,P'_n,P'(n))$ but then $A(n)=A'(n)$ since $F'$ was uniformized.

$x=x'$: also follows from the fact that $F'$ was uniformized.

\noindent{\sc Existence.} Now with transfinite induction we construct for each
$\alpha \in \omega_1$ a
$(z_{\alpha},A_{\alpha},P_{\alpha},p_{\alpha},x_{\alpha}) \in F'''$ with the
required properties. 

Let us formulate the inductive hypothesis: let $\alpha<\omega_1$ be an ordinal and suppose that for every $\beta<\alpha$ we have $(z_{\beta},A_{\beta},P_{\beta},p_{\beta},x_{\beta}) \in F'''$ such that for every $\beta<\alpha$ we have $\{A_{\beta}(n):n \in \omega\}=\{x_{\gamma}:\gamma<\beta\}$.

We will construct $(z_{\alpha},A_{\alpha},P_{\alpha},p_{\alpha},x_{\alpha}) \in F'''$ satisfying the previous hypothesis.

$z_{\alpha}$: using the $1^{st}$ point of Lemma \ref{fi} there exists a unique $z_{\alpha}$ such that $<_{z_{\alpha}} \cong \alpha$ and $(\exists x \in \om)\phi(z_{\alpha},x)$. 

$p_{\alpha}$: let $p_{\alpha}$ be the $\alpha^{th}$ real with respect to $<_L$.

$A_{\alpha}, P_{\alpha}$: The order-preserving bijection between $<_{z_{\alpha}}$ and $\alpha$ yields enumerations $\{x_{\beta}:\beta<\alpha\}$ and $\{p_{\beta}:\beta<\alpha\}$, let $A_{\alpha}(n)$ be the $n^{th}$ element of the first set's enumeration and define $P_{\alpha}(n)$ similarly. 

By the definition of $A_{\alpha}$ we have that $\{A_{\alpha}(n):n \in \omega\}=\{x_{\beta}:\beta<\alpha\}$.

We will prove that there exists an $x_{\alpha}\in \om$ such that
$(z_{\alpha},A_{\alpha},P_{\alpha},p_{\alpha},x_{\alpha}) \in F'''$. By the
properties of $F$ for every $(A,p)$ there exist cofinaly many (in the hyperdegrees)
$x$ such that $(A,p,x)\in F$ , so this also holds for
$(A_{\alpha},p_{\alpha})$. From Lemma \ref{l2} we have
that if the section $F'_{(z_{\alpha},A_{\alpha},P_{\alpha},p_{\alpha})}$ is
non-empty then it is cofinal in the hyperdegrees. 

Now we show that it is non-empty. $L \models $"$P_{\alpha}$ is an enumeration of the
first $\alpha$ reals given by $<_{z_{\alpha}}$ and $p_{\alpha}$ is the
$\alpha^{th}$ real" so by absoluteness arguments it holds in
$L_{\omega^x_1}$ is $\omega^x_1$ is high enough. Let us choose a real $x$ such that $x \in F_{A_{\alpha},p_{\alpha}} \cap \mathcal{S}$, $L_{\omega^x_1} \models $"$P_{\alpha}$ is an enumeration of the
first $\alpha$ reals given by $<_{z_{\alpha}}$ and $p_{\alpha}$ is the $\alpha^{th}$ real" and $\phi(z_{\alpha},x)$. Such an $x$ exists by the $2^{nd}$ point of Lemma \ref{fi} and by the fact that $F_{(A,p)} \cap \mathcal{S}$ is cofinal in the hyperdegrees. Clearly $(z_{\alpha},A_{\alpha},P_{\alpha},p_{\alpha},x) \in F'$. 

Thus there exists an $x_{\alpha}$ such that $(z_{\alpha},A_{\alpha},P_{\alpha},p_{\alpha},x_{\alpha}) \in F''$.

What remains to show is that $(z_{\alpha},A_{\alpha},P_{\alpha},p_{\alpha},x_{\alpha}) \in F'''$:

 From $(z_{\alpha},A_{\alpha},P_{\alpha},p_{\alpha},x_{\alpha}) \in F'$ follows that $\phi(z_{\alpha},x_{\alpha})$. First notice that by the $4^{th}$ 
point of Lemma \ref{fi} 
$\phi(z_{\alpha},x_{\alpha})$ implies the existence of $y_n$-s and $g_n$-s satisfying properties $2(a)$ and $2(b)$ from the definition of $F'''$. 

To see that 
$2(c)$ also holds for
$(z_{\alpha},A_{\alpha},P_{\alpha},p_{\alpha},x_{\alpha})$, fix a natural number
$n$. We know that $\phi(y_n,x_{\alpha})$ holds thus there exists a $\beta<\alpha$ such
that $<_{y_n} \cong \beta$. For all $\beta < \alpha$ the formula $\phi(z_{\beta},x_{\beta})$ holds (by
inductive hypothesis $(z_{\beta},A_{\beta},P_{\beta},p_{\beta},x_{\beta}) \in
F''' \subset F'$ and use the $1^{st}$ point of the definition of $F'$). Let us set $A_n=A_{\alpha} \circ g_n$ and $P_n=P_{\alpha} \circ g_n$. 

We will prove that 
\[(y_n,A_n,P_n,P_{\alpha}(n),A_{\alpha}(n))=(z_{\beta},A_{\beta},P_{\beta},p_{
\beta},x_{\beta}) \in F''. \]

By the $1^{st}$ property of $\phi$ the equality
$y_n=z_{\beta}$ holds. 

Now using the inductive hypothesis we have that $\{A_{\beta}(m):m \in
\omega\}=\{x_{\gamma}:\gamma<\beta\}$. The latter set clearly equals $\{A_n(m):m \in \omega\}$. $A_{\beta}$ and $A_n$ are the enumerations of
the same set of reals given by $<_{z_{\beta}}=<_{y_n}$, hence $A_n=A_{\beta}$. 

Similarly, since $P_{\beta}$ and $P_n$ are the enumerations of the same set
(namely the $\beta$ long initial segment of the reals with respect to $<_L$, see
the Existence part of the proof and Remark \ref{meg}).
Finally, $A_{\alpha}(n)$ and $P_{\alpha}(n)$ are defined as $x_{\beta}$ and the
$\beta^{th}$ real, respectively. 

This finishes the proof of the statement that $2(c)$ also holds for
$(z_{\alpha},A_{\alpha},P_{\alpha},p_{\alpha},x_{\alpha})$ and hence the
proof of the existence. \qed

We have already seen that $X$ is a $\Pi^1_1(t)$ set. Now we check that it is compatible with $F$.
By the previous proposition, for every $\alpha<\omega_1$ there exists a unique
element $(z_{\alpha},A_{\alpha},P_{\alpha},p_{\alpha},x_{\alpha}) \in F'''$ such
that $<_{z_{\alpha}} \cong \alpha$. This gives us the enumerations
$X=\{x_{\alpha}:\alpha<\omega_1\}$ and $\{p_{\alpha}:\alpha<\omega_1\}$. Now by
$3^{rd}$ point of the definition of $F'$ we have that if
$(z_{\alpha},A_{\alpha},P_{\alpha},p_{\alpha},x_{\alpha}) \in F''' \subset F'$
then $L_{\omega^{x_{\alpha}}_1} \models $'$p_{\alpha}$ is the $\alpha^{th}$ real
with respect to $<_L$' and by absoluteness the same holds in $L$. Thus we
obtain that $\om=\{p_{\alpha}:\alpha<\omega_1\}$.  Fix an $\alpha <\omega_1$. 
By the second claim of Proposition \ref{aaa} it is clear that $A_{\alpha}$ is an
enumeration of $\{x_{\beta}:\beta<\alpha\}$. Furthermore, 
$(z_{\alpha},A_{\alpha},P_{\alpha},p_{\alpha},x_{\alpha}) \in F''' \subset F'$ thus by the $2^{nd}$ point of the definition of $F'$ we have that
 $x_{\alpha} \in F_{(A_{\alpha},p_{\alpha})}$, so we can conclude that $X$ is compatible with $F$. \qed
\section{Generalizations and remarks}
Now we will prove the following theorem. 
\thm{}{($V=L$) Let $B$ be a Borel subset of an arbitrary Polish space,
$|B|>\aleph_0$.
Suppose that $F \subset (\om)^{\leq \omega} \times B \times \om$ is a
coanalytic set and for all $p \in B$, $A \in (\om)^{
\leq \omega}$ the section $F_{(A,p)}$ is cofinal in the hyperdegrees. Then there
exists a coanalytic set $X \subset \om$  that is compatible with $F$.
}

\proof A classical result states that for every uncountable Borel subset
$B$ of
a Polish space there exists a map $\Psi: \om \to B$ that is a Borel
isomorphism. 

Suppose that $F$ is a set as above. Let us define $G \subset (\om)^{\leq \omega} \times \om \times \om$ as follows 
\[(A,q,x) \in G \iff (A,\Psi(q),x) \in F. \]
 Clearly, $G$ is a coanalytic set thus there exists a $t \in \om$ so that $G \in
\Pi^1_1(t)$. Of course, each section $G_{(A,q)}$ is cofinal in the hyperdegrees.
The direct application of the Main Theorem yields a $\Pi^1_1(t)$ (therefore
coanalytic) set $X \subset \om$ that is compatible with $G$. From the
compatibility
we obtain the enumeration $\om=\{q_{\alpha}:\alpha<\omega_1\}$. But then
$\{\Psi(q_{\alpha}):\alpha<\omega_1\}$ is an enumeration of $B$ and clearly, $X$
is compatible with $F$ using this enumeration. \qed

\bigskip

We can derive an obvious but useful consequence of the previous theorem using that
$x \leq_T y$ implies $x \leq_h y$ and omitting the relativization.
\thm{turing}{($V=L$) Let $P$ be an uncountable Borel subset of a Polish space.
Suppose that $F \subset (\om)^{\leq \omega} \times P
\times \om$ is a coanalytic set and for all $p \in \om$, $A \in (\om)^{\leq
\omega}$ the section $F_{(A,p)}$ is cofinal in the Turing degrees.
Then there exists a coanalytic set $X$ that is compatible with $F$.}

It is also easy to see that in the previous theorem we can replace $\om$ by $\R^n$ or $2^{\omega}$ etc., since there are recursive Borel isomorphisms between these 
spaces. Thus we obtain Theorem \ref{alap}.

With the same
methods one could prove the following strengthening of the Main Theorem:
\thm{}{($V=L$) Let $B$ be a $\Delta^1_1(t)$ subset of $\om$, $|B|>\aleph_0$.
Suppose that $F \subset (\om)^{\leq \omega} \times B \times \om$ is a
$\Pi^1_1(t)$ set and for all $p \in B$, $A \in (\om)^{
\leq \omega}$ the section $F_{(A,p)}$ is cofinal in the hyperdegrees. Then there exists an $X \in \Pi^1_1 (t)$ that is compatible with $F$.
}

\bigskip
Now we will examine the necessity of $(V=L)$. 
\thm{}{If the conclusion of the Main Theorem holds then there exists a
$\Sigma^1_2$ well-ordering of the reals. In particular, every real is
constructible.}
\proof
Fix recursive $\Delta^1_1$ bijections $\Psi_1:\om \to (\om)^{\leq \omega} \times
\om$ and $\Psi_2:\om \to \om \times \om$.

Let us define the set $F \subset (\om)^{\leq \omega} \times \om \times \om$ as follows:

$(A,p,x) \in F \iff (A,p)=\Psi_1(\pi_{1}(\Psi_2(x)) \land (\forall n)(A(n) \not = x),$

where $\pi_1$ is the projection of $\om \times \om$ on the first coordinate. So
basically $x$ is coding the previous choices and the parameter in the "odd
coordinates".
 
$F$ is clearly $\Delta^1_1$. Now for an arbitrary pair $(A,p)$ and $y \in \om$
there exist cofinaly many $x \in \om$ such that
$(A,p)=\Psi_1(\pi_1(\Psi_2(x))$ and $y \leq_h x$,  hence every section
$F_{(A,p)}$ is cofinal in the hyperdegrees. Thus by our hypothesis there exists
a $\Pi^1_1$ set $X=\{x_{\alpha}:\alpha<\omega_1\}$ and an enumeration
$\om=\{p_{\alpha}:\alpha<\omega_1\}$ such that for every $\alpha<\omega_1$ we
have $x_{\alpha} \in F_{(A_{\alpha},p_{\alpha})}$, where $A_{\alpha}$ is an
enumeration of $\{x_{\beta}:\beta<\alpha\}$. 

We will define the well-ordering of $\om$ with the help of the given enumeration
of $X$.  Since every $x_{\alpha}$ codes the appropriate $p_{\alpha}$, we can
order $\om$ by the first appearance of a real $p$. 

Now for $p,q \in \om$ let $(p,q) \in E \iff \exists x,y,A,B$ 

\begin{enumerate}
\item $x,y \in X$, $x \not =y$, $(A,p,x) \in F, (B,q,y) \in F$
\item $(\forall m)(\forall C)((C,p,A(m)) \not \in F \land (C,q,B(m)) \not \in F)$
\item $(\exists n)(x=B(n))$.
\end{enumerate}

Since $F$ is $\Delta^1_1$, we have that $E$ is a $\Sigma^1_2$ relation. 

Fix $p, q \in \om$. There exist minimal ordinals $\alpha, \beta$ such that $p_{\alpha}=p$ and $p_{\beta}=q$. We will prove that $(p,q) \in E \iff \alpha<\beta$. We have for $\alpha$ and $\beta$
that $(A_{\alpha},p_{\alpha},x_{\alpha}) \in F$ and $(A_{\beta},p_{\beta},x_{\beta}) \in F$. 

First, if $\alpha<\beta$ choose $x=x_{\alpha}$, $y=x_{\beta}$, $A=A_{\alpha}$,
$B=A_{\beta}$. Then $1$ is obvious (by the definition of $F$
we have that $x_{\alpha} \not =x_{\beta}$ if $\alpha<\beta$) and $A_{\beta}$ is
an enumeration of
$\{x_{\gamma}:\gamma<\beta\}$ so $3$ also holds. Suppose that $2$
fails for $p$: there exists a pair $m, C$ such that $(C,p,A(m)) \in F$
(the other case is similar). Then $A(m)=x_{\gamma}$ for some $\gamma<\alpha$ and
$(C,p)=(A_{\gamma},p_{\gamma})$. This would contradict the minimality of
$\alpha$, and similarly for $\beta$.

For the other direction suppose that $(p,q) \in E$ and take $x,y,A,B$ witnessing
this fact. Clearly, $x=x_{\alpha'}$ for some $\alpha'$ so
$(A_{\alpha'},p_{\alpha'})=(A,p)$ and similarly $(A_{\beta'},p_{\beta'})=(B,q)$.
Using $2$ we get the minimality of $\alpha'$ and $\beta'$ so they must be equal
to $\alpha$ and $\beta$. 

Suppose that $\alpha \geq \beta$, then of course $\alpha>\beta$. By $3$ we have
that there exists an $n \in \omega$ such that
\[A_{\beta}(n)=A_{\beta'}(n)=B(n)=x=x_{\alpha'}=x_{\alpha}.\] 
By the assumption $\{x_{\gamma} :\gamma<\beta\} \subsetneq
\{x_{\gamma}:\gamma<\alpha\}$. We have that
\[\{A_{\beta}(m): m \in \omega\}=\{x_{\gamma}:\gamma<\beta\} \subset
\{A_{\alpha}(m): \in \omega\} \]then $A_{\alpha}(m)=x_{\alpha}$ for some $m \in
\omega$. But this is a contradiction, since $(\forall n)(A(n) \not = x)$ for
every $(A,p,x) \in F$. Thus $\alpha<\beta$.  

So we obtain that $E$ is a $\Sigma^1_2$ well-ordering. The second claim follows from Mansfield's theorem, see \cite[Theorem 25.39]{jech}. 
\qed

\bigskip
Next we show that the definability assumption on our ``selection algorithm'' $F$
cannot be dropped in the Main Theorem.

\ex{}{(CH) There exists a family $\{A_{\alpha}:\alpha<\omega_1\} \subset
[\om]^{\leq \aleph_0}$ such that if for a set $X$ there exists an enumeration
$X=\{x_{\alpha}:\alpha<\omega_1\}$ so that $(\forall \alpha
<\omega_1)(x_{\alpha} \not \in A_{\alpha})$ then $X$ is not coanalytic.}

\proof  Fix an enumeration of the reals $\{y_{\alpha}:\alpha<\omega_1\}$. We
will define $A_{\alpha}$ by recursion. Suppose that we are ready for
$\beta<\alpha$ and let us choose $A_{\alpha} \in [\om]^{\leq \aleph_0}$ such
that for every uncountable $P \in \bigcup_{\beta \leq \alpha}
\Pi^1_1(y_{\beta})$ we have $|P \cap (A_{\alpha} \setminus
\bigcup_{\beta<\alpha} A_{\beta})| \geq 2$ and $\bigcup_{\beta<\alpha} A_{\beta}
\subset A_{\alpha}$ and $y_{\alpha} \in A_{\alpha}$. Since
$|\bigcup_{\beta<\alpha} A_{\beta}| \leq \aleph_0$ and $\bigcup_{\beta \leq
\alpha}
\Pi^1_1(y_{\beta})$ is countable, there exists such an
$A_{\alpha}$. 

Now suppose that $X=\{x_{\alpha}:\alpha<\omega_1\}$ is coanalytic and for every
$\alpha$ we have $x_{\alpha} \not \in A_{\alpha}$. Clearly, $\bigcup_{\alpha}
A_{\alpha}=\om$, thus $X$ must be uncountable. Since $X$ is coanalytic, we have
that there exist an $\alpha_0$ such that $X \in \Pi^1_1(y_{\alpha_0})$. Thus for
every $\alpha \geq \alpha_0$ by the construction of $A_{\alpha}$'s $|X \cap
(A_{\alpha} \setminus \bigcup_{\beta<\alpha} A_{\beta})| \geq 2$. Now consider
the map $\phi$ that assigns to each $\alpha \geq \alpha_0$ the minimal index
$\phi(\alpha)$ such that $x_{\phi(\alpha)} \in A_{\alpha+1} \setminus
A_{\alpha}$. There are at least two distinct elements of $X$ in $A_{\alpha+1}
\setminus A_{\alpha}$ and $x_{\gamma} \not \in A_{\alpha+1}$ for
$\gamma>\alpha$ (the constructed family is increasing), hence
$\phi(\alpha)<\alpha$. Moreover, $\phi$ is clearly injective. Therefore, we have
that $\phi$ is a regressive function whose domain is a co-countable subset of
$\omega_1$. This contradicts Fodor's lemma. 
\qed

\rem{}{The same holds for any projective class.}

\bigskip
Now we will prove a general technical theorem which implies the existence of $\Pi^1_1$ Hamel basis, but could be used to prove the existence of $\Pi^1_1$ 
$n$-point sets, analogous versions for circles, etc. The situation in the following definition is that we have a relation $R(x,y)$ on finite subsets of the reals 
that intuitively means that $x$ is "stronger" than $y$ in some sense (e.g. in
case of Hamel basis all elements of $y$ are linearly generated by $x$, in case
of 
two-point sets 
all lines that intersect $y$ in at least two points intersect $x$ in at least two points etc.). Our goal is to find an $R$-independent set (all the relations are trivial) that is "stronger" than all the finite subsets of the reals. $H^R_B$ will be the set of finite sets that can be added to $B$ preserving it's independence.
\defi{}{Let $R$ be a binary relation on the finite subsets of $\R^n$. 
\begin{itemize}
\item We say that a set $X \subset \R^n$ is $R$-independent if for all $x,y \in [X]^{<\omega}$ $R(x,y) \Rightarrow y \subset x$. 
\item Fix a $k \in \omega$, if for every $y \in [\R^n]^k$ there exists an element $x \in [X]^{<\omega}$ such that $R(x,y)$ then we say that $X$ is a $k$-generator set for $R$. 
\item If $B$ is an $R$-independent set let us use the notation $H^R_B=\{x\in [\R^n]^{<\omega}:x \cup B$ is $R$-independent$\}$.
\end{itemize}}
We use parameters $n$ and $k$ even though they will not be needed for the proof of the Hamel basis case. 
\defi{}{We will use the following notation: $x \equiv_h y \iff (x \leq_h y \land y \leq_h x)$.}

The extra difficulty in the construction of a Hamel basis is that in a step we have to put more than one real into our set, so we have to deal with finite sequences. Moreover, to use our method one have to choose reals which are high enough in $\leq_h$. Thus our strategy is to select $\leq_h$ equivalent reals in every step of the procedure.
\defi{}{Let us denote by $\mathcal{E}$ the set \[\{x \in [\R^n]^{<\omega}:(\forall x_1,x_2 \in x)(x_1 \equiv_h x_2)\}.\]}
\thm{teknik}{($V=L$) Let $t \in \R$ and $n,k \in \omega$ be arbitrary. Suppose
that $R \subset [\R^n]^{<\omega} \times [\R^n]^{<\omega}$ is a $\Delta^1_1(t)$ 
relation that satisfies the property (*): \begin{center}for every countable $B \subset \R^n$ the set $\mathcal{E} \cap H^R_B $ is cofinal in the hyperdegrees 
and if for $y \in [\R^n]^k$ there is no $z \in [B]^{<\omega}$ such that $R(z,y)$ then $\{x: R(x,y)\}\cap \mathcal{E} \cap H^R_B$ is cofinal in the hyperdegrees. \end{center}
Then there exists an uncountable $\Pi^1_1(t)$, $R$-independent set that is a $k$-generator for $R$. }

\proof Let us define the set $F \subset ([\R^n]^{<\omega})^{\leq \omega} \times \R \times [\R^n]^{<\omega}$ and fix a recursive Borel isomorphism $\Phi: \R \to 
[\R^n]^k$. \\
$(A,p,x) \in F \iff$\\
EITHER the conjunction of the following clauses holds
\begin{enumerate}
 \item $\bigcup ran(A)$ is $R$-independent
 \item $(\forall z \in ran(A))(\lnot R(z,\Phi(p)))$
 \item $R(x,\Phi(p))$ holds and $x\in \mathcal{E} \cap H^R_{\bigcup ran(A)} $ 
\end{enumerate}
OR $1 \land \lnot 2$ holds and $x \in \mathcal{E} \cap H^R_{\bigcup ran(A)}$\\
OR $\lnot 1$.

Since $A$ is countable and the relation $\equiv_h$ is $\Pi^1_1$, we get that $F$ is $\Pi^1_1(t)$. By property (*) every section $F_{(A,p)}$ is cofinal in the 
hyperdegrees (if $\lnot 1$ then this is obvious and the cases when $1 \land \lnot 2$ or $1 \land 2$ holds are exactly described by property (*))
so we can apply Theorem \ref{main}. This gives us a $\Pi^1_1(t)$ set $Y \subset [\R^n]^{<\omega}$ such that $\bigcup ran(Y)$ is $R$-independent and for every $y \in 
[\R^n]^k$ 
there exists an $x \in Y$ such that $R(x,y)$ thus $\bigcup ran(Y)$ is a $k$-generator for $R$. Moreover $ran(Y) \subset \mathcal{E}$. Hence it suffices  to prove that 
$X=\bigcup ran(Y)$ is a $\Pi^1_1(t)$ set. But using that for every $x \in Y$ the elements of $x$ are equivalent in hyperdegrees we get
\[ a \in X \iff (\exists l \in \omega)(\exists a_1,\dots a_l \leq_h a)(\{a,a_1\dots a_l\} \in ran(Y)).\]
Applying Theorem \ref{sck} we can verify that $X \in \Pi^1_1(t)$.
\qed

\cor{}{($V=L$) There exists a $\Pi^1_1$ Hamel basis.}
\proof
Let us define the relation $R \subset [\R]^{<\omega} \times [\R]^{<\omega}$. $R(x,y) \iff (y \subset \langle x \rangle_{\mathbb{Q}})$ i. e. every element of $y$ is in the linear subspace generated by the elements of $x$ over the rationals. Notice that $R$ is $\Delta^1_1$. In the terminology of the previous theorem $X$ is a Hamel basis if it is $R$-independent and $1$-generator for $R$. So we just have to check whether property (*) holds.

First if $B$ is a countable linearly independent subset of the reals then for
all but countably many finite sets $a \in [\R]^{<\omega}$ we have $a \in H^R_B$.
Therefore obviously $H^R_B$ is cofinal in the hyperdegrees. So the first part of
(*) holds.

Now fix an element $y \in \R$, a countable $B \subset \R$ such that there is no
$z \in [B]^{<\omega}$ such that $R(z,\{y\})$. We will prove that for every $s
\in \R$ there exists a pair $w_1,w_2 \in \R$ satisfying $y=w_1+w_2$, $w_1
\equiv_h w_2$, $B \cup \{w_1,w_2\}$ linearly independent and $s \leq_h w_1,
w_2$. This fact indeed implies that the set $\{x: x \in \mathcal{E} \land
R(x,y)\} \cap H^R_B$ is cofinal in the hyperdegrees, so the second part of (*)
also holds.

 Here we repeat Miller's argument. Without loss of generality we can suppose that $y \leq_h s$ and $s$ is not hyperarithmetic in any finite subset of $B \cup \{y\}$ because we can replace $s$ by a more complicated real. We can choose $w_1$ and $w_2$ such that $s$ is coded in $w_1$'s odd and $w_2$'s even digits so that $w_1+w_2=y$. Then $s \leq_h w_1,w_2$ hence $y \leq_h w_1, w_2$. But then $y=w_1+w_2$ implies $w_1 \equiv_h w_2$. If $w_1 \in \langle B, w_2 \rangle_{\mathbb{Q}}$ then $y \in \langle B, w_2 \rangle_{\mathbb{Q}} \setminus \langle B \rangle_{\mathbb{Q}}$ and then $w_2 \in \langle B,y \rangle_{\mathbb{Q}}$ but this would imply that $s$ is hyperarithmetic in a finite subset of $B \cup \{ y\}$ which is a contradiction. Thus $w_1$ and $w_2$ are the appropriate reals.

Thus property (*) holds indeed, and the direct application of Theorem \ref{teknik} hence produces a $\Pi^1_1$  Hamel basis. 
\qed

\bigskip
Finally we will prove another variant of our theorem, considering the case where the choice at step $\alpha$ does not depend on the previous choices.

\thm{UDC}{($V=L$) Let be $t \in \R$ and suppose that $G \subset \R^n \times \R$
is a $\Delta^1_1(t)$ set and for every countable $A \subset \R$ the complement 
of 
the set $\cup_{p \in A} G_p$ is cofinal in the hyperdegrees. Then there exists an uncountable $\Pi^1_1(t)$ set $X \subset \R^n$ that intersects every $G_p$ in a countable set.}
\proof
Using Theorem \ref{spg} there exists a $\Sigma_1$ formula $\theta$ such that \[a \in G^c \iff L_{\omega^{(a,t)}_1}[a,t] \models \theta(a,t).\] 

Now let us define the set $H$ as follows:
\[(x,p) \in H \iff x \in \mathcal{S} \land p,t \leq_h x \land L_{\omega^x_1} \models (\forall p' \leq_L p)(\theta((x,p'),t)).\]

$H$ is a $\Pi^1_1(t)$ set, for this just repeat the usual argument, that is, $x \in \mathcal{S} \land p,t \leq_h x$ implies that 
$L_{\omega^x_1}=L_{\omega^{((x,p),t)}_1}[((x,p),t)]$ and use Theorems \ref{spg}, \ref{a}, \ref{b} and Lemma \ref{TEK}. Observe that for a real $p$  
\[H_p=(\bigcap_{p' \leq_L p} G^c_{p'}) \cap \mathcal{S} \cap \{z: p,t \leq_h z\}. \]
Thus the theorem's conditions imply that for every real $p$ the section $H_p$ is cofinal in the hyperdegrees.

Define $F \subset (\R^n)^{\leq \omega} \times \R \times \R^n$:
$(A,p,x) \in F \iff (x,p) \in H \land x \not \in A$. Obviously for every $(A,p)$ the section $F_{(A,p)}$ is cofinal in the hyperdegrees and $F$ is $\Pi^1_1(t)$. 
Our Main Theorem provides an uncountable $\Pi^1_1(t)$ set $X \subset \R^n$  and enumerations $X=\{x_{\alpha}:\alpha<\omega_1\}$, 
$\R=\{p_{\alpha}:\alpha<\omega_1\}$ and an enumeration $A_{\alpha}$ (in type $\leq \omega$) of $\{x_{\beta}:\beta<\alpha\}$ such that $x_{\alpha} \in 
F_{(A_{\alpha},p_{\alpha})}=H_{p_{\alpha}} \setminus \{x_{\beta}:\beta<\alpha\}$. Suppose that there exists a $p \in \R$ for which $|X \cap G_p|>\aleph_0$. Then $p_{\beta}>_L p$ if $\beta$ is high enough, since only countably many $p_{\alpha}$'s are $<_L$ less then $p$. But if $p_{\beta}>_L p$ then $x_{\beta} \in G^c_p$. \qed

Now Theorem \ref{UDC0} is a trivial consequence of Theorem \ref{UDC}.

\section{Applications}

 Theorem \ref{alap} can be applied in various situations. Let us remark here that one can obtain $\Pi^1_1$ sets instead of coanalytic ones by just repeating the proofs and using Theorem \ref{main} in all the theorems of this section. We will prove the simpler (boldface) versions for the sake of transparency.
\thm{}{($V=L$) There exists a coanalytic MAD family.}
\proof First fix a recursive partition $B=\{B_i:i \in \omega\}$ of $\omega$ to infinite sets.
Define $F \subset (\mathcal{P}(\omega))^{\leq \omega} \times \mathcal{P}(\omega)
\times \mathcal{P}(\omega)$ as follows: $(A,p,x) \in F \iff$\\EITHER the
conjunction of the following clauses holds
\begin{enumerate}
\item $ran(A) \cup B$ contains pairwise almost disjoint elements
\item $p$ is almost disjoint form the elements of $ran(A) \cup B$
\item $p \subset x$ and $x$ is almost disjoint form the elements of $ran(A) \cup B$
\end{enumerate}
OR $1 \land \lnot2$ holds and $x$ is almost disjoint form the elements of $ran(A) \cup B$\\
OR $ \lnot 1.$

Clearly, $F$ is Borel. What we have to prove is that for all pairs $(A,
p)$ the section $F_{(A,p)}$ is cofinal in the Turing degrees. 

Suppose that $1$ and $2$ hold, let $u \in \mathcal{P}(\omega)$ be an arbitrary real. Choose $x'=p \cup \bigcup_{i \in \omega} F_i$, 
where $F_i \subset B_i$ are finite and if $i>j$ then $A(j) \cap F_i=\emptyset$ and \[|(p \cup F_i) \cap B_i|\equiv 1 \mod 2 \iff u(i)=1.\]
For every $i$ there exist such an $F_i$, since the $B_i$'s are disjoint and infinite, and $ran(A)\cup B$ contains pairwise almost disjoint sets. Then $x'$ satisfies $3$ and $u \leq_T x'$. 

Now in the case when $1 \land \lnot 2$ holds our job is easier: e. g. we can
repeat the previous argument omitting $p$.

Finally, if $\lnot 1$ is true then $F_{(A,p)}=\mathcal{P}(\omega)$.

 Notice that Theorem \ref{alap} was stated in the form that the set of the parameters is $\R$ but we can easily replace it by $\mathcal{P}(\omega)$ using a recursive Borel isomorphism. 

So we can apply Theorem \ref{alap} and we get a coanalytic set $X=\{x_{\alpha}:\alpha \in \omega_1\}$ such that $X$ is compatible with $F$.  
It is obvious by transfinite induction that the elements of $X$ are pairwise almost disjoint. It is also clear that $X \cup B$ is maximal since for every real $p$ there exists an $\alpha<\omega_1$ such that $p_{\alpha}=p$. Thus there exists an element of $X$ that is not almost disjoint from $p$.\qed

\thm{}{($V=L$) There exists a coanalytic two-point set.}

\proof For each real $p \in \mathbb{R}$ fix a line $l_p$ such that it is the line defined by the equation $((p)_1)x+((p)_2)y=(p)_3$, where $(p)_1,(p)_2$ and $(p)_3$ are the reals made of every $3k^{th}, 3k+1^{th}$ and $3k+2^{th}$ digit of $p$. $l_p$ can be empty, however every line appears at least two times. Let us define $F \subset (\mathbb{R}^{2})^{\leq \omega} \times \R \times \R^2$ by
$(A,p,x) \in F \iff $\\
EITHER the conjunction of the following clauses holds
\begin{enumerate}
\item there are no $3$ collinear points in $ran(A)$
\item $|ran(A) \cap l_p|<2$ and $l_p \not = \emptyset$
\item $x \in l_p \setminus ran(A)$ and $x$ is not collinear with any two distinct points of $ran(A)$
\end{enumerate}
OR $1 \land \lnot 2$ holds and $x$ is not collinear with two distinct points of $ran(A)$\\
OR $\lnot 1$.

Now $F$ is clearly Borel. What we have to check is that for all $(A,p)$ the section $F_{(A,p)}$ is cofinal in the Turing degrees. Fix a pair $(A,p)$. If $1 \land 
2$ holds then the section is equal to $l_p$ minus a countable set. Every line is cofinal in the Turing degrees, because we can choose one of the coordinates 
arbitrarily. Now notice that if $H$ is a set which is cofinal in the Turing degrees and $H'$ is countable the $H \setminus H'$ is still cofinal: to see this let $u$ be an arbitrary real and let $s$ be such that $(\forall s' \in H')(s' \not \geq_T s)$ then there exist $r \in H$ such that $s, u \leq_T r$ and clearly $r \not \in H'$. So we have that if $1 \land 2$ holds then $F_{(A,p)}$ is cofinal in the Turing degrees.

If $1 \land \lnot 2$ holds then we just have to choose an arbitrary point that
is not 
collinear with any two distinct points of $A$. The 
case when $1$ is false is obvious.  

Thus by Theorem \ref{alap} we get an uncountable coanalytic set $X=\{x_{\alpha}:\alpha<\omega_1\} \subset \R^2$. One can easily verify that $X$ cannot contain three collinear points. Moreover, since every line $l_p$ appears at least twice, $|l_p \cap X|=2$.
\qed

Similar statements can be formulated for $n$-point sets, circles, appropriate algebraic curves etc., the above method works in these cases. 
\subsection{Curves in the plane}
Now we will consider the following question: What can we say about a set in the plane which intersects every "nice" curve in a countable set?
Let us call a continuously differentiable $\R \to \R^2$ function a $C^1$ curve.

\defi{c1}{We say that a set $H \subset \R^2$ is $C^1$-small if the
intersection of $H$ with the range of every $C^1$ curve is a countable set.}

In \cite{kun} the authors proved that assuming Martin's axiom and the Semi-Open
Coloring Axiom if $H$ is $C^1$-small then $|H| \leq \aleph_0$. Moreover, they
showed in ZFC that no perfect set is $C^1$-small. Thus no uncountable analytic
set is $C^1$-small. On the other hand, the following proposition holds. 

\prop{}{(CH) There exists an uncountable $C^1$-small set.}
\proof We will prove later that the union of the range of countably many $C^1$ curves cannot cover the plane. This implies the statement by an easy transfinite induction.  \qed

Thus it is interesting whether an uncountable $C^1$-small subset can be coanalytic. We will 
apply Theorem \ref{UDC0}.
\thm{}{($V=L$) There exists an uncountable $C^1$-small coanalytic set.}
\proof
First we have to prove that there exists a Borel set $G \subset \R^2 \times \R$ such that if $\gamma$ is a $C^1$ curve then there exists a $p \in \R$ such that $G_p=ran(\gamma)$.

One can easily prove that the set $B$ of $C^1$ curves as a subset of $C(\R,\R^2)$ is a Borel set (see e.g.  \cite[23. D]{cdst}). The set $\{((x,y),\gamma):(x,y) \in 
ran(\gamma)\} \subset \R^2 \times C(\R,\R^2)$ is clearly closed. So $(\R^2 \times B) \cap \{((x,y),\gamma):(x,y) \in ran(\gamma)\}$ is also a Borel set. Furthermore, there exists a Borel isomorphism $\phi: \R \to B$ since these two are standard Borel spaces of cardinality $\mathfrak{c}$ and we can apply the isomorphism theorem. Now we can define $G \subset \R^2 \times \R$: $((x,y),p) \in G \iff ((x,y),\phi(p)) \in (\R^2 \times B) \cap \{((x,y),\gamma):(x,y) \in ran(\gamma)\}$ which is a Borel set and for every $\gamma \in C^1$ there exists a $p \in \R$ such that $G_p=ran(\gamma)$.

To apply Theorem \ref{UDC0} we have to check that if we have countably many $C^1$ curves $\{ \gamma_i:i \in \omega\}$ then the complement of the union of their ranges is cofinal in the Turing degrees. For this it is enough that there exists a line $l$ such that \[|l \cap \bigcup (\{ ran(\gamma_i):i\in \omega\})| \leq \aleph_0.\] Let us concentrate solely on the horizontal lines. For a curve $\gamma_i$ take let $f_i(x)=\pi_y(\gamma_i(x))$, i. e. the composition with the projection on the vertical axis. $f_i$ is $C^1$ function, thus by Sard's lemma the set $H_i=\{y \in \R: (\exists x)(f'_i(x)=0 \land f_i(x)=y\}$ has Lebesgue measure zero. Let $b \in \R \setminus (\cup H_i)$. Then the line $\{(x,b):x\in \R\}$ intersects every curve $\gamma_i$ in countably many points, since otherwise it would be an image of a critical value.

Finally, the application of Theorem \ref{UDC0} produces an uncountable $C^1$-small coanalytic set.

\qed

\subsection{Problems}

In Theorem \ref{alap} the set of the parameters is a Borel set and this was used
in the proof numerous times. 

\prb{}{Does Theorem \ref{alap} hold if we only assume that $B$ is coanalytic?}

As a partial converse we have proved that the conclusion of the Main Theorem
implies that every real constructible. It is natural to ask whether the
converse also holds.

\prb{}{Does the conclusion of Theorem \ref{alap} hold if every real is
constructible?}

One of the weaknesses of the method is that the constructed set $X$ is
a subset of $\mathcal{S}$. It is known (see e. g. \cite{kechris}) that
$\mathcal{S}$ is the largest thin (not containing a perfect
subset) $\Pi^1_1$ set. Thus non of the constructed sets contain a perfect
subset. In the case of $C^1$-small sets this cannot be expected, but how
about the other constructions?

\prb{}{Is it consistent that there exists a $\Pi^1_1$ Hamel basis (two-point
set, MAD family) that contains a perfect subset?}

\textbf{Acknowledgement} I am very greatful to my supervisor, M\'arton Elekes, for his patience and the help what he provided during the
writing of this paper.

\bigskip

\noindent
\textsc{Institute of Mathematics, E\"otv\"os Lor\'and University, 
P\'azm\'any P\'eter s. 1/c, Budapest 1117, Hungary}

\noindent
\textit{Email address}: \verb+vidnyanszkyz@gmail.com+

\noindent
{\tt www.cs.elte.hu/\hbox{$\sim$}vidnyanz}

\end{document}